\documentclass[journal]{IEEEtran}
\ifCLASSINFOpdf
\else
\fi
\usepackage{cite}
\usepackage{amssymb}
\usepackage{amsmath}
\usepackage{graphicx}
\newtheorem{assumption}{Assumption}
\newtheorem{lemma}{Lemma}
\newtheorem{theorem}{Theorem}
\newtheorem{definition}{Definition}
\newtheorem{remark}{Remark}

\usepackage{epstopdf}
\usepackage{color}
\usepackage{cases}

\usepackage{subfigure}

\usepackage{bm}

\begin{document}
\title{Aug-PDG: Linear Convergence of Convex Optimization with Inequality Constraints}

\author{
Min Meng and Xiuxian Li% <-this % stops a space
\thanks{M. Meng (mengmin@tongji.edu.cn) and X. Li (xli@tongji.edu.cn) are with the Department of Control Science and Engineering, College of Electronics and Information Engineering, Tongji University, Shanghai, 201804, China and the Institute for Advanced Study, Tongji University, Shanghai, 200092, China.
}% <-this % stops a space
%\thanks{Corresponding author: }
%\thanks{This work was partially supported by .}% <-this % stops a space
}

\maketitle

\begin{abstract}
This paper investigates the convex optimization problem with general convex inequality constraints. To cope with this problem, a discrete-time algorithm, called augmented primal-dual gradient algorithm (Aug-PDG), is studied and analyzed. It is shown that Aug-PDG can converge semi-globally to the optimizer at a linear rate under some mild assumptions, such as the quadratic gradient growth condition for the objective function, which is strictly weaker than strong convexity. To our best knowledge, this paper is the first to establish a linear convergence for the studied problem in the discrete-time setting, where an explicit bound is provided for the stepsize. Finally, a numerical example is presented to illustrate the efficacy of the theoretical finding.
\end{abstract}


% \fi
\IEEEpeerreviewmaketitle

\section{Introduction}
This paper deals with the constrained optimization problem formulated as follows:
\begin{align}
&\min\limits_{x\in\mathbb{R}^n}~~ f(x)\nonumber\\
&\text{~~s.t.} ~ ~g(x)\leq 0,\label{equ1}
\end{align} where the objective function $f:\mathbb{R}^n\rightarrow\mathbb{R}$ and $g(x)=(g_1(x),g_2(x),\ldots,g_m(x))^T$ with $g_i:\mathbb{R}^n\rightarrow\mathbb{R}$ being convex and continuously differentiable. By resorting to the (or augmented) Lagrangian $L(x,\lambda)$ of problem (\ref{equ1}), the corresponding (or augmented) primal-dual gradient algorithm (PDG) (or Aug-PDG) can be designed as
 \begin{align}
 x_{k+1}&=x_k-\alpha\nabla_xL(x_k,\lambda_k),\nonumber\\
 \lambda_{k+1}&=[\lambda_k+\alpha\nabla_{\lambda}L(x_k,\lambda_k)]_+,\label{equ2}
 \end{align} where $\alpha$ is a positive stepsize and $[\cdot]_+$ denotes the projection operator onto the nonnegative orthants component-wisely. It is known that (\ref{equ2}) can find a saddle point of the Lagrangian $L(x,\lambda)$, and thus it has been extensively studied to solve the constrained optimization problem \cite{ruszczynski2011nonlinear}.

Optimization has wide applications in power systems \cite{chiang2007layering,zhao2014design}, wireless communication \cite{chen2012convergence}, game theory \cite{gharesifard2013distributed,pavel2020distributed}, to name just a few. To date, there is a large body of literature on theoretical analysis of asymptotic convergence of various algorithms, including primal-dual gradient-based algorithms, for tackling the optimization problem under different settings \cite{polyak1970iterative,golshtein1974generalized,nedic2009subgradient,feijer2010stability,cherukuri2016asymptotic,hong2017prox,lei2017doubly,li2019distributed,li2020distri}.

In recent decades, researchers have focused on the exponential/linear convergence of primal-dual gradient-based algorithms.
It is well-known that when the objective function is strongly convex and smooth, the gradient decent algorithm for unconstrained convex optimization can achieve global exponential convergence in continuous-time and global linear convergence in discrete-time. In the context of constrained optimization with equality constraints $Ax= b$ or affine inequality constraints  $Ax\leq b$, PDG is proved to converge globally and exponentially in continuous-time setup \cite{niederlander2016exponentially}. A proximal gradient
flow was proposed in \cite{dhingra2019proximal}, which can be applied to resolve convex optimization problems with affine inequality constraints and has global exponential convergence when $A$ has full row rank. Local exponential convergence of the primal-dual gradient dynamics can be established with the help of spectral bounds of saddle matrices \cite{shen2010eigenvalue}. Recently, the authors in \cite{qu2019exponential} proved that the Aug-PDGD in continuous-time for optimization with affine equality and inequality constraints achieves global exponential convergence, and the global linear converge of primal-dual gradient optimization (PDGO) in discrete-time was discussed in \cite{su2019contraction} by contraction theory.
It should be noted that the aforementioned works focus on unconstrained optimization or constrained optimization with affine equality and/or affine inequality constraints. For the case with nonlinear inequality constraints, the asymptotic/sublinear convergence has been extensively studied such as in \cite{nedic2011random}. However, the linear convergence for the optimization with nonlinear inequality constraints is seldom investigated in the literature. One exception is the recent work \cite{tang2019semi}, where the authors established a semi-global exponential convergence of Aug-PDGD in the sense that the convergence rate depends on the distance from the initial point to the optimal point.

%A new algorithm relying on the subgradient projection was discussed for the optimization problem with set and nonlinear inequality constraints in \cite{nedic2011random}, while only sublinear convergence was ensured.
% More recently, \cite{tang2019semi} extended the method in \cite{qu2019exponential} to nonlinear inequality constraints and obtained

However, \cite{tang2019semi} concentrates on the continuous-time dynamics. As discrete-time algorithms are easily implemented in practical applications, in this paper, the discrete-time algorithm is addressed for the optimization problem with nonlinear inequality constraints. Theoretical analysis based on a quadratic Lyapunov function that has non-zero off-diagonal terms is presented to show that the Aug-PDG achieves semi-global linear convergence, where an explicit bound is established for the stepsize. The numerical results suggest that the Aug-PDG indeed has different linear convergence rates for different initial points.

The rest of this paper is organized as follows. Section \ref{sec2} introduces preliminaries on optimization with nonlinear equality constraints. The main result on the semi-global linear convergence of Aug-PDGA, along with its proof, is presented in Section \ref{sec3}. Section \ref{sec4} provides a numerical example to illustrate the feasibility of the obtained result. Section \ref{sec5} makes a brief conclusion.

{\em Notations.}  Let $\mathbb{R}^m$, $\mathbb{R}_+^m$ and $\mathbb{R}^{m\times n}$ be the sets of $m$-dimensional real column vectors, $m$-dimensional nonnegative column vectors and $m\times n$ real matrices, respectively. Define $[x]_+$ to be the component-wise projection of a vector $x\in\mathbb{R}^m$ onto $\mathbb{R}_+^m$. The symbol $x\geq0$ for any vector $x\in\mathbb{R}^m$ means that each entry of $x$ is nonnegative. For an integer $n>0$, denote $[n]:=\{1,2,\ldots,n\}$. $I_n$ is the identity matrix of dimension $n$. ${\bf 1}_n$ (resp. ${\bf 0}_n$) represents an $n$-dimensional vector with all of its elements being 1 (resp. 0). For a vector or matrix $A$, $A^\top$ denotes the transpose of $A$ and $A_{\mathcal I}$ is a matrix composed of the rows of $A$ with the indices in ${\mathcal I}$.
%$\rho(A)$ represents the spectral radius of $A$ and $det(A)$ is the determinant of $A$.
For real symmetric matrices $P$ and $Q$, $P\succ(\succeq, \succ,\preceq)~Q$ means that $P-Q$ is positive (positive semi-, negative, negative semi-) definite, while for two vectors/matrices $w,v$ of the same dimension, $w\leq v$ means that each entry of $w$ is no greater than the corresponding one of $v$. ${\rm diag}\{a_1,a_2,\ldots,a_n\}$ represents a diagonal matrix with $a_i$, $i\in[n]$, on its diagonal.
%$A\otimes B$ denotes the Kronecker product of matrices $A$ and $B$. For a vector $v$, we use ${\rm diag}(v)$ to represent the diagonal matrix with the vector $v$ on its diagonal.

\section{Preliminaries}\label{sec2}
Consider problem (\ref{equ1}). An augmented Lagrangian associated with problem (\ref{equ1}) is introduced as \cite{bertsekas2014constrained}
\begin{align}\label{equ4}
L(x,\lambda):=f(x)+U_\rho(x,\lambda),
\end{align} where $x\in\mathbb{R}^n$, $\lambda=(\lambda_1,\lambda_2,\ldots,\lambda_m)^{\top}\in\mathbb{R}^m$, $\rho>0$ is the penalty parameter, and
\begin{align}\label{equ5}
U_\rho(x,\lambda):=\sum\limits_{i=1}^m\frac{[\rho g_i(x)+\lambda_i]_+^2-\lambda_i^2}{2\rho}.
\end{align} It can be verified that $U_\rho(x,\lambda)$ is convex in $x$ and concave in $\lambda$, and $U_\rho(x,\lambda)$ is continuously differentiable, i.e.,
\begin{align}
\nabla_xU_{\rho}(x,\lambda)&=\sum\limits_{i=1}^m[\rho g_i(x)+\lambda_i]_+\nabla g_i(x),\label{equ6}\\
\nabla_{\lambda}U_{\rho}(x,\lambda)&=\sum\limits_{i=1}^m\frac{[\rho g_i(x)+\lambda_i]_+-\lambda_i}{\rho}e_i,\label{equ7}
\end{align} where $e_i$ is an $n$-dimensional vector with the $i$th entry being 1 and others 0. Then the augmented primal-dual gradient algorithm (Aug-PDG) can be explicitly written as
\begin{subequations}\label{equ8}
\begin{align}
x_{k+1}&=x_k-\alpha\nabla_xL(x_k,\lambda_k)\nonumber\\
&=x_k-\alpha\nabla f(x_k)-\alpha\sum\limits_{i=1}^m[\rho g_i(x_k)+\lambda_{i,k}]_+\nabla g_i(x_k),\\
\lambda_{k+1}&=\lambda_k+\alpha\nabla_{\lambda}L(x_k,\lambda_k)\nonumber\\
&=\lambda_k+\alpha\sum\limits_{i=1}^m\frac{[\rho g_i(x_k)+\lambda_{i,k}]_+-\lambda_{i,k}}{\rho}e_i,\label{equ8b}
\end{align}
\end{subequations} where $\alpha\in(0,\rho]$ is the stepsize to be specified. Here, the initial conditions are arbitrarily chosen as $x_0\in\mathbb{R}^n$ and $\lambda_0\geq 0$.

To proceed, the following results are vital for solving the constrained optimization problem.
\begin{lemma}\label{lem2}
For Aug-PDG (\ref{equ8}), if $\lambda_0\geq0$, then $\lambda_k\geq0$ for all $k\geq0$.
\end{lemma}

{\em Proof.} This result can be proved by mathematical induction. First note that $\lambda_0\geq 0$. Assume now that $\lambda_k\geq0$ for some $k\geq0$, then by (\ref{equ8b}), one has that
\begin{align*}
\lambda_{k+1}&=(1-\alpha/\rho)\lambda_k+\alpha\sum\limits_{i=1}^m\frac{[\rho g_i(x_k)+\lambda_{i,k}]_+}{\rho}e_i\\
&\geq(1-\alpha/\rho)\lambda_k\\
&\geq 0,
\end{align*} where the first inequality is obtained based on the definition of the notation $[\cdot]_+$ and the second inequality is derived following $0<\alpha\leq\rho$ and inductive assumption $\lambda_k\geq0$. The proof is thus completed. \hfill$\blacksquare$
\begin{lemma}
A primal-dual pair $(x^*,\lambda^*)$ is an equilibrium point of the Aug-PDG (\ref{equ8}) if and only if $(x^*,\lambda^*)$ is a Karush-Kuhn-Tucker (KKT) point of (\ref{equ1}).
\end{lemma}

{\em Proof.} If a primal-dual pair $(x^*,\lambda^*)$ is an equilibrium point of the Aug-PDG (\ref{equ8}), that is,
\begin{align*}
x^*&=x^*-\alpha\nabla_xL(x^*,\lambda^*),\\
\lambda^*&=\lambda^*+\alpha\nabla_{\lambda}L(x^*,\lambda^*),
\end{align*} then $\nabla_xL(x^*,\lambda^*)=0$ and $\nabla_\lambda L(x^*,\lambda^*)=0$. For $\nabla_\lambda L(x^*,\lambda^*)=0$, an equivalent condition is that for any $i\in[m]$,
\begin{align}\label{equ9}
\lambda_i^*=[\rho g_i(x^*)+\lambda_i^*]_+,
\end{align} which implies $\lambda_i^*\geq0$, $g_i(x^*)\leq 0$, and $\lambda_i^* g_i(x^*)=0$. For $\nabla_xL(x^*,\lambda^*)=0$, one can equivalently obtain that $\nabla f(x^*)+\sum_{i=1}^m[\rho g_i(x^*)+\lambda_i^*]_+\nabla g_i(x^*)=\nabla f(x^*)+\sum_{i=1}^m\lambda_i^*\nabla g_i(x^*)=0$. Thus, it can be claimed that the primal-dual pair $(x^*,\lambda^*)$ is a KKT point.

Conversely, if $(x^*,\lambda^*)$ is a KKT point of (\ref{equ1}), then $\nabla f(x^*)+\sum_{i=1}^m\lambda_i^*\nabla g_i(x^*)=0$, $\lambda_i^*g_i(x^*)=0$, $\lambda^*\geq0$ and $g_i(x^*)\leq 0$. Via a simple computation, $\nabla_xL(x^*,\lambda^*)=0$ and $\nabla_{\lambda}L(x^*,\lambda^*)=0$, which implies that $(x^*,\lambda^*)$ is an equilibrium point of the Aug-PDG (\ref{equ8}). \hfill$\blacksquare$
%
%\begin{remark}
%The augmented Lagrangian used in (\ref{equ4}) is different from the standard Lagrangian in \cite{feijer2010stability} by employing $U(x,\lambda)$ in (\ref{equ5}) instead of $\lambda^{\top}g(x)$. In addition, the standard PDGA involves a discontinuous projection onto the nonnegative orthant, which may create difficulties both in theoretic analysis and numerical simulations. It can been seen from Lemma \ref{lem2} that $\lambda_k$ can be automatically guaranteed to be nonnegative provided that the initial value $\lambda_0$ is nonnegative, thus avoiding performing projection and discontinuity issues caused by the projection step.
%\end{remark}
\section{Main Results}\label{sec3}
In this section, the main result on the linear convergence of Aug-PDG is presented.
\subsection{Convergence Results}
The following assumptions are essential for deriving the main result.
\begin{assumption}\label{assump1}
The problem (\ref{equ1}) has a unique feasible solution $x^*$, and at $x^*$, the linear independence constraint qualification (LICQ) holds at $x^*$, i.e., $\{\nabla g_i(x^*)\mid i\in{\mathcal I}\}$ is linearly independent, where ${\mathcal I}:=\{i\in[m]\mid g_i(x^*)=0\}$ is the so-called active set at $x^*$.
\end{assumption}

Under Assumption \ref{assump1}, the optimal Lagrangian multiplier $\lambda^*$ is also unique \cite{wachsmuth2013licq}. Denote by $J$ the Jacobian of $g(x)$ at $x^*$ and $J_{\mathcal I}$ the matrix composed of the rows of $J$ with the indices in ${\mathcal I}$. LICQ in Assumption \ref{assump1} also implies that $J_{\mathcal I}J_{\mathcal I}^{\top}\succ0$ \cite{tang2019semi}. Define
\begin{align}
\kappa:=\lambda_{\min}(J_{\mathcal I}J_{\mathcal I}^{\top})>0
\end{align} to be the smallest eigenvalue of $J_{\mathcal{I}}J_{\mathcal{I}}^\top$.

\begin{assumption}\label{assump2}
The objective function $f(x)$ has a quadratic gradient growth with parameter $\mu>0$ over $\mathbb{R}^n$, i.e.,
\begin{align}
(\nabla f(x)-\nabla f(x^*))^{\top}(x-x^*)\geq\mu\|x-x^*\|^2,~\forall x\in \mathbb{R}^n.
\end{align}
\end{assumption}

The concept of quadratic gradient growth was introduced in \cite{necoara2019linear}, which is a relaxation of strong convexity condition for guaranteeing linear convergence of gradient-based optimization algorithms. In fact, the class of functions having quadratic gradient growth include the strongly convex functions as a proper subset and some functions with quadratic gradient growth are even not convex.

\begin{assumption}\label{assump3}
The objective function $f$ is $l$-smooth over $\mathbb{R}^n$, i.e., $\|\nabla f(x)-\nabla f(y)\|\leq l\|x-y\|$ for any $x,y\in \mathbb{R}^n$. For any $i\in[m]$, $g_i(x)$ is $L_{gi}$-smooth and has bounded gradient, i.e., $\|\nabla g_i(x)-\nabla g_i(y)\|\leq L_{gi}\|x-y\|$ and $\|\nabla g_i(x)\|\leq B_{gi}$ for some $L_{gi},B_{gi}>0$ and any $x,y\in \mathbb{R}^n$.
\end{assumption}

Denote ${\mathcal I}^c:=[m]\backslash{\mathcal I}$, $L_g:=\sqrt{\sum_{i=1}^mL_{gi}^2}$ and $B_g:=\sqrt{\sum_{i=1}^mB_{gi}^2}$. Under Assumption \ref{assump3}, one can obtain that
\begin{align}
\|J\|&\leq B_g,\label{equ11}\\
\|g(x)-g(y)\|&=\sqrt{\sum_{i=1}^m(g_i(x)-g_i(y))^2}\nonumber\\
&\leq\sqrt{\sum_{i=1}^mB_{gi}^2\|x-y\|^2}\nonumber\\
&=B_g\|x-y\|. \label{equ10}
\end{align}

Denote $d_0:=\sqrt{\|x_0-x^*\|^2+\|\lambda_0-\lambda^*\|^2}$. Before giving the main result of this paper, it is convenient to list the following concept similar to that in continuous-time setting \cite{sastry2013nonlinear}.
\begin{definition}
Consider the dynamics $z(t+1)=\phi(z(t))$ with initial point $z(0)=z_0$. Assume that $z_e$ is an equilibrium point satisfying $z_e=\phi(z_e)$. $z_e$ is said to be a semi-global linear stable point if for any $h>0$, there exist $c>0$ and $0<\gamma<1$ such that for any $z_0$ satisfying $\|z_0-z_e\|\leq h$,
\begin{align*}
\|z(t)-z_e\|\leq c\gamma^t\|z_0-z_e\|,~\forall t\geq0.
\end{align*} $z_e$ is said to be a global linear stable point if $c$ and $\gamma$ do not depend on $h$.
\end{definition}

Then the main result is presented as follows.
\begin{theorem}\label{thm1}
Under Assumptions \ref{assump1}--\ref{assump3}, if the stepsize $\alpha>0$ is chosen such that
\begin{align}
\alpha<\min\left\{1,\rho,\frac{2\mu}{b_1+2a_4\delta},\frac{\kappa\delta}{2b_2+4a_5\delta},\frac{1-\pi^*}{2\rho(b_2+2a_5\delta)}\right\},
\end{align} where $\delta>0$ satisfies
\begin{align}
\delta<\min\left\{\frac{\mu}{2a_3},\frac{1-\pi^*}{2\rho(\kappa+8B_g^2+L_g^2(1-\pi^*))},B_g^{-1}\right\},\label{equ14}
 \end{align}
 $\pi^*:=[\rho \max_{i\in{\mathcal I}^c}\{g_i(x^*)\}/(\sqrt{C}d_0)+1]_+^2$, $b_1:=a_1+2B_g^2$, $b_2:=a_2+\frac{2}{\rho^2}$, $a_1:=2l+4\theta_1^2$, $a_2:=4B_g^2$, $a_3:=2B_g^2l^2/\kappa+2B_g^2\theta_1^2/\kappa+2B_g^2/(\kappa\rho^2)+{\kappa}B_g^2\rho^2/4$, $a_4:=B_g^2l^2/2+B_g^2\theta_1^2+2B_g^2$, $a_5:=B_g^2+2/\rho^2$, and $\theta_1:=\rho B_g^2+L_g\|\lambda^*\|$, then the sequences $\{x_k\}$ and $\{\lambda_k\}$ generated by Aug-PDG (\ref{equ8}) for the constrained optimization (\ref{equ1}) semi-globally converge to the optimal point of the optimization problem (\ref{equ1}) at a linear (or exponential) rate. Specifically, it holds that
\begin{align}\label{equ15}
\|x_{k}-x^*\|^2+\|\lambda_{k}-\lambda^*\|^2\leq C(1-\gamma)^kd_0^2,
\end{align} where $0<\gamma<1$ satisfies
\begin{align}\label{e16}
\gamma\leq\min\{c_1,c_2,c_3\}
\end{align} with $c_1:=\mu\alpha-a_3\delta\alpha-b_1\alpha^2/2-a_4\delta\alpha^2$,
$c_2:=\kappa\delta\alpha/4-b_2\alpha^2/2-a_5\delta\alpha^2$,
$c_3:=\frac{\alpha}{2\rho}(1-\pi^*)-(\delta\alpha\kappa+b_2\alpha^2+2a_5\delta\alpha^2)/2-4\alpha\delta B_g^2$, and $C:=\lambda_{\max}(Q_\delta)/\lambda_{\min}(Q_\delta)\geq1$ with $Q_\delta:=\left[\begin{array}{cc}I_n&\delta J^{\top}\\ \delta J & I_m\end{array}\right]$.
\end{theorem}

Proof. The proof is postponed to the next subsection. \hfill$\blacksquare$
\begin{remark}
The selection of parameters $\alpha$ and $\delta$ ensures that $c_1,c_2,c_3$ are positive, and then $\gamma>0$ can be guaranteed. From Theorem \ref{thm1}, one can see that the convergence rate is related to $\pi^*=[\rho \max\{g_i(x^*)\}/(\sqrt{C}d_0)+1]_+^2$ and decrease to 0 as $d_0$ goes to infinity.  The decreasing rate also changes as $(x_k,\lambda_k)$ approaches the optimal point. Specifically, the decreasing rates are small at the beginning and then become large when $(x_k,\lambda_k)$ goes to the optimal point. Therefore, Theorem \ref{thm1} does not guarantee the existence of a global linear convergence rate, and consequently only semi-global linear stability can be ensured.
\end{remark}
\begin{remark}
To our best knowledge, this paper is the first to investigate the linear convergence for problem (\ref{equ1}) in discrete-time setup. Compared with the most related literature \cite{tang2019semi}, where a continuous-time algorithm, called Aug-PDGD, was studied with a semi-global exponential convergence, a discrete-time algorithm Aug-PDG is analyzed here with a semi-global linear convergence. Although discrete-time algorithms may be obtained by discretizing the continuous-time Aug-PDGD using such as explicit Euler method, it is unclear how to select the sampling stepsize to guarantee the convergence especially in the sense of semi-global convergence. In comparison, an explicit bound on the stepsize $\alpha$ is established here in Theorem \ref{thm1}.
\end{remark}
\subsection{Proof of Theorem 1}
To prove Theorem \ref{thm1}, an intermediate result is needed as follows.
\begin{lemma}\cite{qu2019exponential}\label{lem3}
For any $y,y^*\in\mathbb{R}$, there exists $\xi\in[0,1]$ such that $[y]_+-[y^*]_+=\xi(y-y^*)$. Specifically, $\xi$ can be chosen as $\xi=\frac{[y]_+-[y^*]_+}{y-y^*}$ if $y\neq y^*$ and $\xi=0$ if $y=y^*$.
\end{lemma}

%{\em Proof.} The proof is obvious, which is thus omitted. \hfill$\blacksquare$
Then the proof of Theorem \ref{thm1} is presented as follows.

Define
\begin{align}
{V}_{\delta,k}&=\left[\begin{array}{cc}{x}_k-x^*\\ \lambda_k-\lambda^*\end{array}\right]^{\top}Q_\delta\left[\begin{array}{c}{x}_k-x^*\\ \lambda_k-\lambda^*\end{array}\right],
\end{align} where
\begin{align}
Q_{\delta}=\left[\begin{array}{cc}I_n&\delta J^{\top}\\ \delta J & I_m\end{array}\right].
\end{align} As $\delta<B_g^{-1}$ from (\ref{equ14}), one has $I_m-\delta^2JJ^{\top}\succeq(1-\delta^2B_g^2)I_m\succ0$, which implies  $Q_\delta\succ0$ by Schur complement. Then, in the following, we discuss the bound of
\begin{align*}
{V}_{\delta,k+1}=&\|{x}_{k+1}-x^*\|^2+\|\lambda_{k+1}-\lambda^*\|^2\\
&+2\delta({x}_{k+1}-x^*)^{\top}J^{\top}(\lambda_{k+1}-\lambda^*).
\end{align*}

Note that $(x^*,\lambda^*)$ is the KKT point of (\ref{equ1}), that is,
\begin{align}
\nabla_x L(x^*,\lambda^*)&=0,\\
\nabla_\lambda L(x^*,\lambda^*)&=0.
\end{align}
%
%Based on (\ref{equ4})--(\ref{equ7}), one has that
%\begin{align}
%&\nabla_xL(x_k,\lambda_k)\nonumber\\
%&=\nabla_xL(x_k,\lambda_k)-\nabla_xL(x^*,\lambda^*)\nonumber\\
%&=\nabla f(x_k)-\nabla f(x^*)+\nabla_xU(x_k,\lambda_k)-\nabla_xU(x^*,\lambda^*)
%\end{align}

By iterations in (\ref{equ8}), one has that
\begin{align}
&\|{x}_{k+1}-x^*\|^2\nonumber\\
&=\|x_k-\alpha\nabla_xL(x_k,\lambda_k)-x^*+\alpha\nabla_xL(x^*,\lambda^*)\|^2\nonumber\\
&=\|x_k-x^*\|^2+\alpha^2\|\nabla_xL(x_k,\lambda_k)-\nabla_x L(x^*,\lambda^*)\|^2\nonumber\\
&~~~-2\alpha(\nabla_xL(x_k,\lambda_k)-\nabla_xL(x^*,\lambda^*))^{\top}(x_k-x^*).\label{equ18}
\end{align}  By $\nabla_xL(x_k,\lambda_k)=\nabla f(x_k)+\nabla_xU(x_k,\lambda_k)=\nabla f(x_k)+\sum_{i=1}^m[\rho g_i(x_k)+\lambda_{i,k}]_+\nabla g_i(x_k)$, for the second term on the right side of (\ref{equ18}), one has that
\begin{align}
&\alpha^2\|\nabla_xL(x_k,\lambda_k)-\nabla_x L(x^*,\lambda^*)\|^2\nonumber\\
&\leq2\alpha^2\|\nabla f(x_k)-\nabla f(x^*)\|^2\nonumber\\
&~~~+2\alpha^2\|\nabla_xU(x_k,\lambda_k)-\nabla_xU(x^*,\lambda^*)\|^2.\label{equ19}
\end{align} Note that
\begin{align*}
&\nabla_xU(x_k,\lambda_k)-\nabla_xU(x^*,\lambda^*)\\
&=\sum_{i=1}^m[\rho g_i(x_k)+\lambda_{i,k}]_+\nabla g_i(x_k)\\
&~~~-\sum_{i=1}^m[\rho g_i(x^*)+\lambda_{i}^*]_+\nabla g_i(x^*)\\
&=\sum_{i=1}^m([\rho g_i(x_k)+\lambda_{i,k}]_+-[\rho g_i(x^*)+\lambda_{i}^*]_+)\nabla g_i(x_k)\\
&~~~~+\sum_{i=1}^m[\rho g_i(x^*)+\lambda_{i}^*]_+(\nabla g_i(x_k)-\nabla g_i(x^*)),
\end{align*} then
\begin{align}
&\|\nabla_xU(x_k,\lambda_k)-\nabla_xU(x^*,\lambda^*)\|\nonumber\\
&\leq\sum_{i=1}^m|[\rho g_i(x_k)+\lambda_{i,k}]_+-[\rho g_i(x^*)+\lambda_{i}^*]_+|\cdot\|\nabla g_i(x_k)\|\nonumber\\
&~~~~+\sum_{i=1}^m[\rho g_i(x^*)+\lambda_{i}^*]_+\|\nabla g_i(x_k)-\nabla g_i(x^*)\|.\label{equ20}
\end{align} Define
\begin{align}
\xi_{i,k}=\frac{[\rho g_i(x_k)+\lambda_{i,k}]_+-[\rho g_i(x^*)+\lambda^*]_+}{(\rho g_i(x_k)+\lambda_{i,k})-(\rho g_i(x^*)+\lambda^*)}\label{equ21}
\end{align} if $(\rho g_i(x_k)+\lambda_{i,k})-(\rho g_i(x^*)+\lambda^*)\neq0$, and
\begin{align}
\xi_{i,k}=0\label{equ22}
\end{align} if $(\rho g_i(x_k)+\lambda_{i,k})-(\rho g_i(x^*)+\lambda^*)=0$. Then it can be obtained from Lemma \ref{lem3} that
\begin{align}
&[\rho g_i(x_k)+\lambda_{i,k}]_+-[\rho g_i(x^*)+\lambda_{i}^*]_+\nonumber\\
&=\rho\xi_{i,k}(g_i(x_k)-g_i(x^*))+\xi_{i,k}(\lambda_{i,k}-\lambda_i^*).\label{equ23}
\end{align} Substituting (\ref{equ23}) into (\ref{equ20}) yields that
\begin{align}
&\|\nabla_xU(x_k,\lambda_k)-\nabla_xU(x^*,\lambda^*)\|\nonumber\\
&\leq\sum_{i=1}^m|\rho\xi_{i,k}(g_i(x_k)-g_i(x^*))+\xi_{i,k}(\lambda_{i,k}-\lambda_i^*)|\cdot\|\nabla g_i(x_k)\|\nonumber\\
&~~~+\sum_{i=1}^m[\rho g_i(x^*)+\lambda_{i}^*]_+\|\nabla g_i(x_k)-\nabla g_i(x^*)\|\nonumber\\
&\leq \sum_{i=1}^m(\rho B_{gi}\|x_k-x^*\|+|\lambda_{i,k}-\lambda_i^*|)B_{gi}\nonumber\\
&~~~+\sum_{i=1}^m\lambda_i^*L_{gi}\|x_k-x^*\|\nonumber\\
&\leq \theta_1\|x_k-x^*\|+B_g\|\lambda_k-\lambda^*\|,\label{equ24}
\end{align} where Assumption \ref{assump3} has been applied to get the second inequality, and the third inequality is derived by $\sum_{i=1}^na_ib_i\leq\sqrt{\sum_{i=1}^na_i^2}\sqrt{\sum_{i=1}^nb_i^2}$ for any $a_i,b_i\in\mathbb{R}$, $\theta_1=\rho B_g^2+L_g\|\lambda^*\|$ and $B_g=\sqrt{\sum_{i=1}^mB_{gi}^2}$, $L_g=\sqrt{\sum_{i=1}^mL_{gi}^2}$.

By (\ref{equ24}), (\ref{equ19}) and Assumption \ref{assump3}, one has that
\begin{align}
&\alpha^2\|\nabla_xL(x_k,\lambda_k)-\nabla_x L(x^*,\lambda^*)\|^2\nonumber\\
&\leq2\alpha^2l^2\|x_k-x^*\|^2\nonumber\\&~~~+2\alpha^2(\theta_1\|x_k-x^*\|+B_g\|\lambda_k-\lambda^*\|)^2\nonumber\\
&\leq2\alpha^2l^2\|x_k-x^*\|^2\nonumber\\&~~~+4\alpha^2(\theta_1^2\|x_k-x^*\|^2+B_g^2\|\lambda_k-\lambda^*\|^2)\nonumber\\
&=a_1\alpha^2\|x_k-x^*\|^2+a_2\alpha^2\|\lambda_k-\lambda^*\|^2,\label{equ25}
\end{align} where $a_1=2l^2+4\theta_1^2$ and $a_2=4B_g^2$.

For the third term on the right side of (\ref{equ18}),
\begin{align}
&-2\alpha(\nabla_xL(x_k,\lambda_k)-\nabla_xL(x^*,\lambda^*))^{\top}(x_k-x^*)\nonumber\\
&=-2\alpha(\nabla f(x_k)-\nabla f(x^*))^{\top}(x_k-x^*)\nonumber\\
&~~~-2\alpha(\nabla_xU(x_k,\lambda_k)-\nabla_xU(x^*,\lambda^*))^{\top}(x_k-x^*)\nonumber\\
&\leq-2\mu\alpha\|x_k-x^*\|^2+2\alpha(U(x^*,\lambda_k)-U(x_k,\lambda_k))\nonumber\\
&~~~+2\alpha(U(x_k,\lambda^*)-U(x^*,\lambda^*)),\label{equ26}
\end{align} where the inequality is derived based on Assumption \ref{assump2} and the convexity of $U(x,\lambda)$ at $x$, i.e.,
\begin{align}
U(x,\lambda)-U(x',\lambda)\geq(\nabla_xU(x',\lambda))^{\top}(x-x')
\end{align} for any $x,x'\in \mathbb{R}^n$.

Plugging (\ref{equ25}) and (\ref{equ26}) into (\ref{equ18}), one can obtain that
\begin{align}
&\|{x}_{k+1}-x^*\|^2\nonumber\\
&\leq(1-2\mu\alpha+a_1\alpha^2)\|x_k-x^*\|^2+a_2\alpha^2\|\lambda_k-\lambda^*\|^2\nonumber\\
&~~~+2\alpha(U(x^*,\lambda_k)-U(x_k,\lambda_k))\nonumber\\&~~~+2\alpha(U(x_k,\lambda^*)-U(x^*,\lambda^*)).\label{equ29}
\end{align}

For $\|\lambda_{k+1}-\lambda^*\|^2$, by iteration (\ref{equ8b}), one has that
\begin{align}
&\|\lambda_{k+1}-\lambda^*\|^2\nonumber\\
&=\|\lambda_k+\alpha\nabla_\lambda U(x_k,\lambda_k)-\lambda^*\|^2\nonumber\\
&=\|\lambda_k-\lambda^*\|^2+\alpha^2\|\nabla_\lambda U(x_k,\lambda_k)\|^2\nonumber\\
&~~~+2\alpha\nabla_\lambda U(x_k,\lambda_k))^{\top}(\lambda_k-\lambda^*).\label{equ30}
\end{align} Recalling $\nabla_{\lambda}L(x^*,\lambda^*)=\nabla_{\lambda}U(x^*,\lambda^*)=0$ and the notation of $\xi_{i,k}$ in (\ref{equ21}), (\ref{equ22}), it can be obtained that
{\small\begin{align}
&\|\nabla_\lambda U(x_k,\lambda_k)\|^2\nonumber\\
&=\|\nabla_\lambda U(x_k,\lambda_k)-\nabla_\lambda U(x^*,\lambda^*)\|^2\nonumber\\
&=\left\|\sum_{i=1}^m\frac{[\rho g_i(x_k)+\lambda_{i,k}]_+-\lambda_{i,k}-[\rho g_i(x^*)+\lambda_{i}^*]_++\lambda_{i}^*}{\rho}e_i\right\|^2\nonumber\\
&=\left\|\sum_{i=1}^m[\xi_{i,k}(g_i(x_{k})-g_i(x^*))+\frac{1}{\rho}(\xi_{i,k}-1)(\lambda_{i,k}-\lambda_i^*)]e_i\right\|^2\nonumber\\
&=\left\|\Xi_k(g(x_k)-g(x^*))+\frac{1}{\rho}(\Xi_k-I_{m})(\lambda_k-\lambda^*)\right\|^2\nonumber\\
&\leq 2B_g^2\|x_k-x^*\|^2+\frac{2}{\rho^2}\|\lambda_k-\lambda^*\|^2,\label{equ31}
\end{align}}where $\Xi_k=\text{diag}\{\xi_{1,k},\ldots,\xi_{m,k}\}$, the inequality is obtained based on (\ref{equ10}) and $\|\Xi_k\|\leq 1$, $\|\Xi_k-I_m\|\leq1$ for $\xi_{i,k}\in[0,1]$, $i\in[m]$.

In consideration that $U(x,\lambda)$ is concave at $\lambda$, one has
\begin{align}
(\nabla_\lambda U(x_k,\lambda_k))^{\top}(\lambda_k-\lambda^*)\leq U(x_k,\lambda_k)-U(x_k,\lambda^*).\label{equ32}
\end{align}

By (\ref{equ31}) and (\ref{equ32}), it can be derived from (\ref{equ30}) that
\begin{align}
&\|\lambda_{k+1}-\lambda^*\|^2\nonumber\\
&\leq(1+\frac{2\alpha^2}{\rho^2})\|\lambda_k-\lambda^*\|^2+2B_g^2\alpha^2\|x_k-x^*\|^2\nonumber\\
&~~~+2\alpha(U(x_k,\lambda_k)-U(x_k,\lambda^*)).\label{equ33}
\end{align}

On the other hand, it is easy to verify that
\begin{align}
&({x}_{k+1}-x^*)^{\top}J^{\top}(\lambda_{k+1}-\lambda^*)\nonumber\\
&=(x_k-x^*-\alpha\nabla f(x_k)+\alpha\nabla f(x^*)\nonumber\\
&~~~-\alpha\nabla_xU(x_k,\lambda_k)+\alpha\nabla_xU(x^*,\lambda^*))^{\top}J^{\top}\nonumber\\
&~~~(\lambda_k-\lambda^*+\alpha\nabla_\lambda U(x_k,\lambda_k)-\alpha\nabla_\lambda U(x^*,\lambda^*))\nonumber\\
&=(x_k-x^*)^{\top}J^{\top}(\lambda_k-\lambda^*)\nonumber\\
&~~~-\alpha(\nabla f(x_k)-\nabla f(x^*))^{\top}J^{\top}(\lambda_k-\lambda^*)\nonumber\\
&~~~-\alpha(\nabla_xU(x^*,\lambda_k)-\nabla_xU(x^*,\lambda^*))^{\top}J^{\top}(\lambda_k-\lambda^*)\nonumber\\
&~~~-\alpha(\nabla_{x}U(x_k,\lambda_k)-\nabla_xU(x^*,\lambda_k))^{\top}J^{\top}(\lambda_k-\lambda^*)\nonumber\\
&~~~+\alpha({x}_{k+1}-x^*)^{\top}J^{\top}(\nabla_\lambda U(x_k,\lambda_k)-\nabla_\lambda U(x^*,\lambda^*)). \label{equ34}
\end{align}
By (\ref{equ24}), one has that
\begin{align}
\|\nabla_xU(x^*,\lambda_k)-\nabla_xU(x^*,\lambda^*)\|&\leq B_g\|\lambda_k-\lambda^*\|,\label{equ35}\\
\|\nabla_xU(x_k,\lambda_k)-\nabla_xU(x^*,\lambda_k)&\leq\theta_1\|x_k-x^*\|.\label{equ36}
\end{align} Similar to (\ref{equ21})--(\ref{equ22}), define
\begin{align}
\xi_{i,\lambda}=\left\{\begin{array}{ll}\frac{[\rho g_i(x^*)+\lambda_{i,k}]_+-[\rho g_i(x^*)+\lambda^*_i]_+}{\lambda_{i,k}-\lambda^*},&\text{if}~\lambda_{i,k}\neq\lambda_i^*,\\
0,&\text{if}~\lambda_{i,k}=\lambda_i^*,\end{array}\right.
\end{align} then
\begin{align}
&\nabla_xU(x^*,\lambda_k)-\nabla_xU(x^*,\lambda^*)\nonumber\\
&=\sum_{i=1}^m([\rho g_i(x^*)+\lambda_{i,k}]_+-[\rho g_i(x^*)+\lambda_i^*]_+)\nabla g_i(x^*)\nonumber\\
&=\sum_{i=1}^m\xi_{i,\lambda}(\lambda_{i,k}-\lambda_i^*)\nabla g_i(x^*)\nonumber\\
&=J^{\top}\Xi_\lambda(\lambda_k-\lambda^*),\label{equ38}
\end{align} where $\Xi_\lambda:=\text{diag}\{\xi_{1,\lambda},\ldots,\xi_{m,\lambda}\}$. For the last term of (\ref{equ34}), it holds that
{\small\begin{align}
&\alpha({x}_{k+1}-x^*)^{\top}J^{\top}(\nabla_\lambda U(x_k,\lambda_k)-\nabla_\lambda U(x^*,\lambda^*))\nonumber\\
&=\alpha({x}_{k}-x^*)^{\top}J^{\top}(\nabla_\lambda U(x_k,\lambda_k)-\nabla_\lambda U(x^*,\lambda^*))\nonumber\\
&~-\alpha^2(\nabla{f(x_k)}-\nabla{f(x^*)})J^{\top}(\nabla_\lambda U(x_k,\lambda_k)-\nabla_\lambda U(x^*,\lambda^*))\nonumber\\
&~-\alpha^2(\nabla_xU(x_k,\lambda_k)-\nabla_xU(x^*,\lambda^*))^{\top}\nonumber\\&~\times J^{\top}(\nabla_\lambda U(x_k,\lambda_k)-\nabla_\lambda U(x^*,\lambda^*))\nonumber\\
&\leq\frac{2B_g^2\alpha}{\kappa\rho^2}\|x_k-x^*\|^2+\frac{\kappa\rho^2\alpha}{8}\|\nabla_\lambda U(x_k,\lambda_k)-\nabla_\lambda U(x^*,\lambda^*)\|^2\nonumber\\
&~~+\frac{B_g^2\alpha^2}{2}\|\nabla{f(x_k)}-\nabla{f(x^*)}\|^2\nonumber\\
&~~+\frac{B_g^2\alpha^2}{2}\|\nabla_xU(x_k,\lambda_k)-\nabla_xU(x^*,\lambda^*)\|^2\nonumber\\
&~~+\alpha^2\|\nabla_\lambda U(x_k,\lambda_k)-\nabla_\lambda U(x^*,\lambda^*)\|^2,\label{equ43}
\end{align}}where the inequality has been obtained based on $\|J\|\leq B_g$ and $w^{\top}v\leq\frac{\beta}{2}\|w\|^2+\frac{1}{2\beta}\|v\|^2$ for $\beta=\frac{4}{\kappa\rho^2}$ and 1.

Therefore, by (\ref{equ31}), (\ref{equ35})--(\ref{equ43}), one can rewrite (\ref{equ34}) as
\begin{align}
&({x}_{k+1}-x^*)^{\top}J^{\top}(\lambda_{k+1}-\lambda^*)\nonumber\\
&=(x_k-x^*)^{\top}J^{\top}(\lambda_k-\lambda^*)\nonumber\\&-\alpha(\lambda_k-\lambda^*)^{\top}\Xi_{\lambda}JJ^{\top}(\lambda_k-\lambda^*)\nonumber\\
&~~~-\alpha(\nabla f(x_k)-\nabla f(x^*))^{\top}J^{\top}(\lambda_k-\lambda^*)\nonumber\\
&~~~-\alpha(\nabla_xU(x_k,\lambda_k)-\nabla_xU(x^*,\lambda_k))^{\top}J^{\top}(\lambda_k-\lambda^*)\nonumber\\
&~~~+\alpha({x}_{k+1}-x^*)^{\top}J^{\top}(\nabla_\lambda U(x_{k},\lambda_k)-\nabla_\lambda U(x^*,\lambda^*))\nonumber\\
&\leq(x_k-x^*)^{\top}J^{\top}(\lambda_k-\lambda^*)\nonumber\\&~~~-\alpha(\lambda_k-\lambda^*)^{\top}\Xi_{\lambda}JJ^{\top}(\lambda_k-\lambda^*)\nonumber\\
&~~~+\frac{2\|J\|^2\alpha}{\kappa}\|\nabla f(x_k)-\nabla f(x^*)\|^2+\frac{\kappa\alpha}{8}\|\lambda_k-\lambda^*\|^2\nonumber\\
&~~~+\frac{2\|J\|^2\alpha}{\kappa}\|\nabla_xU(x_k,\lambda_k)-\nabla_xU(x^*,\lambda_k)\|^2\nonumber\\
&~~~+\frac{\kappa\alpha}{8}\|\lambda_k-\lambda^*\|^2+\frac{2B_g^2\alpha}{\kappa\rho^2}\|x_k-x^*\|^2\nonumber\\
&~~~+\frac{\kappa\rho^2\alpha}{8}\|\nabla_\lambda U(x_k,\lambda_k)-\nabla_\lambda U(x^*,\lambda^*)\|^2\nonumber\\
&~~~+\frac{B_g^2\alpha^2}{2}\|\nabla{f(x_k)}-\nabla{f(x^*)}\|^2\nonumber\\
&~~~+\frac{B_g^2\alpha^2}{2}\|\nabla_xU(x_k,\lambda_k)-\nabla_xU(x^*,\lambda^*)\|^2\nonumber\\
&~~~+\alpha^2\|\nabla_\lambda U(x_k,\lambda_k)-\nabla_\lambda U(x^*,\lambda^*)\|^2,\nonumber\\
&\leq(x_k-x^*)^{\top}J^{\top}(\lambda_k-\lambda^*)\nonumber\\
&~~~-\alpha(\lambda_k-\lambda^*)^{\top}\Xi_{\lambda}JJ^{\top}(\lambda_k-\lambda^*)\nonumber\\
&~~~+(a_3\alpha+a_4\alpha^2)\|x_{k}-x^*\|^2+(\frac{\kappa\alpha}{2}+a_5\alpha^2)\|\lambda_k-\lambda^*\|^2,\label{equ39}
\end{align} where the last inequality is obtained by a simple computation, along with (\ref{equ24}), (\ref{equ31}) and (\ref{equ36}). %$a_3=2B_g^2l^2/\kappa+2B_g^2\theta_1^2/\kappa+2B_g^2/(\kappa\rho^2)+{\kappa}B_g^2\rho^2/4$, $a_4=B_g^2l^2/2+B_g^2\theta_1^2+2B_g^2$ and $a_5=B_g^2+2/\rho^2$.

Define
\begin{align}\label{equ40}
\tilde{\xi}_{i,k}:=
\left\{\begin{array}{ll}[\rho g_i(x^*)/\lambda_{i,k}+1]_+^2,&i\in{\mathcal I}^c,\lambda_{i,k}\neq0,\\
{[\rho g_i(x^*)/d_0+1]}_+^2,&i\in{\mathcal I}^c,\lambda_{i,k}=0,\\
1,&i\in{\mathcal I},
\end{array}\right.
\end{align} where $\tilde{\Xi}_k:=\text{diag}\{\tilde{\xi}_{1,k},\ldots,\tilde{\xi}_{m,k}\}$. Then
\begin{align}
&U(x^*,\lambda_k)-U(x^*,\lambda^*)\nonumber\\
&=\frac{1}{2\rho}(\lambda_k-\lambda^*)^{\top}(\tilde{\Xi}_k-I_m)(\lambda_k-\lambda^*).\label{equ41}
\end{align}

Combining with (\ref{equ29}), (\ref{equ33}), (\ref{equ39}) and (\ref{equ41}), the bound of ${V}_{\delta,k+1}$ can be derived that
\begin{align}
{V}_{\delta,k+1}&=\|x_{k+1}-x^*\|^2+\|\lambda_{k+1}-\lambda^*\|^2\nonumber\\
&~~~+2\delta({x}_{k+1}-x^*)^{\top}J^{\top}(\lambda_{k+1}-\lambda^*)\nonumber\\
&\leq (1-\gamma)V_{\delta,k}+ \left[\begin{array}{c}x_k-x^*\\ \lambda_k-\lambda^*\end{array}\right]^{\top}Q\left[\begin{array}{c}x_k-x^*\\ \lambda_k-\lambda^*\end{array}\right],\label{equ42}
\end{align}
 where
\begin{align}
Q=\left[
\begin{array}{cc}Q_{1}&\gamma\delta J^{\top}\\ \gamma\delta J&Q_2
\end{array}\right]
\end{align} with $Q_1:=(\gamma-2\mu\alpha+2a_3\delta\alpha+b_1\alpha^2+2a_4\delta\alpha^2)I_n$,
$Q_2:=(\gamma+\delta\alpha\kappa+b_2\alpha^2+2a_5\delta\alpha^2)I_m+\frac{\alpha}{\rho}(\tilde{\Xi}_{k}-I_m)-\alpha\delta(\Xi_{\lambda}JJ^{\top}+JJ^{\top}\Xi_\lambda)$, $b_1=a_1+2B_g^2$, and $b_2=a_2+\frac{2}{\rho^2}$.

If $Q\preceq0$, then ${V}_{\delta,k+1}\leq(1-\gamma)V_{\delta,k}$, indicating that $\lambda_{\min}(Q_\delta)(\|{x}_{k+1}-x^*\|^2+\|{\lambda}_{k+1}-\lambda^*\|^2)\leq\lambda_{\max}(Q_\delta)(1-\gamma)(\|x_k-x^*\|^2+\|\lambda_k-\lambda^*\|^2)$. Therefore, \begin{align}
&\|{x}_{k+1}-x^*\|^2+\|{\lambda}_{k+1}-\lambda^*\|^2\nonumber\\
&\leq\frac{\lambda_{\max}(Q_\delta)}{\lambda_{\min}(Q_\delta)}(1-\gamma)(\|x_k-x^*\|^2+\|\lambda_k-\lambda^*\|^2).\label{equ48}
\end{align} Thus, (\ref{equ15}) holds for $C=\lambda^{-1}_{\min}(Q_\delta)\lambda_{\max}(Q_\delta)$.

Note that
$
\left[\begin{array}{cc}-\gamma I_n&\gamma\delta J^{\top}\\ \gamma\delta J&-\gamma I_m\end{array}\right]\preceq0.
$  Hence, to prove $Q\preceq0$, it suffices to ensure $Q\preceq\left[\begin{array}{cc}-\gamma I_n&\gamma\delta J^{\top}\\ \gamma\delta J&-\gamma I_m\end{array}\right]$, i.e.,
\begin{align}
\left[\begin{array}{cc}Q_1+\gamma I_n&{\bf0}_{n\times m}\\{\bf0}_{m\times n}&Q_2+\gamma I_m\end{array}\right]\preceq0.
\end{align}

%As $\alpha<\frac{\mu}{b_1}$ and $\delta<\frac{a_2\mu}{2}$, then $\mu\alpha-b_1\alpha^2>0$ and $\mu\alpha-2a_2\delta>0$.

By (\ref{e16}), one can obtain that $2\gamma-2\mu\alpha+2a_3\delta\alpha+b_1\alpha^2+2a_4\delta\alpha^2\leq0$, i.e., $Q_1+\gamma I_m\preceq0$.

Next, consider $\Theta:=\frac{\alpha}{\rho}(\tilde{\Xi}_{k}-I_m)-\alpha\delta(\Xi_{\lambda}JJ^{\top}+JJ^{\top}\Xi_\lambda)$ in $Q_2$. If $\Theta+(2\gamma+\delta\alpha\kappa+b_2\alpha^2+2a_5\delta\alpha^2)I_m\preceq0$, then $Q_2+\gamma I_m\preceq0$.

Note that $\xi_{i,\lambda}=\tilde{\xi}_{i,k}=1$ when $i\in{\mathcal I}$. Partition $\Theta$ as
\begin{align}
\Theta=\left[
\begin{array}{cc}
\Theta_{1}&\Theta_3\\ \Theta_3^{\top}&\Theta_2
\end{array}\right]
\end{align} where
\begin{align*}
\Theta_{1}&:=-2\delta\alpha J_{\mathcal I}J_{\mathcal I}^{\top},\\
\Theta_2&:=\frac{\alpha}{\rho}(\tilde{\Xi}_{k,{\mathcal I}^c}-I)-\alpha\delta(\Xi_{\lambda,{\mathcal I}^c}J_{{\mathcal I}^c}J_{{\mathcal I}^c}^{\top}+J_{{\mathcal I}^c}J_{{\mathcal I}^c}^{\top}\Xi_{\lambda,{\mathcal I}^c}),\\
\Theta_3&:=-\delta\alpha J_{\mathcal I}J_{{\mathcal I}^c}^{\top}(I+\Xi_{\lambda,{\mathcal I}^c}).
\end{align*}

Under LICQ in Assumption \ref{assump1}, it can be obtained that $J_IJ_I^{\top}\succeq \kappa I$. Then one can see from (\ref{e16}) that $\Theta_1+(2\gamma+\delta\alpha\kappa+b_2\alpha^2+2a_5\delta\alpha^2)I\preceq-\frac{1}{2}\delta\alpha J_{\mathcal I}{\mathcal J}_{\mathcal I}^{\top}$.

By Lemma 6 in \cite{qu2019exponential}, $L_g^2(I-\Xi_{{\mathcal I}^c})+\Xi_{{\mathcal I}^c}J_{{\mathcal I}^c}J_{{\mathcal I}^c}^{\top}+J_{{\mathcal I}^c}J_{{\mathcal I}^c}^{\top}\Xi_{{\mathcal I}^c}\succeq0$, then since $\tilde{\xi}_{i,k}\leq\xi_{i,\lambda}$ for $i\in{\mathcal I}^c$, one has
\begin{align}
\Theta_2\preceq \frac{\alpha}{\rho}(\tilde{\Xi}_{{\mathcal I}^c}-I)-\delta\alpha L_g^2(\tilde{\Xi}_{{\mathcal I}^c}-I).
\end{align} Denote $\phi=2\gamma+\delta\alpha\kappa+b_2\alpha^2+2a_5\delta\alpha^2$, then one can obtain that
\begin{align}
&\Theta_2+\phi I-\Theta_3^{\top}(\Theta_1+\phi I_m)^{-1}\Theta_3\nonumber\\
&\preceq\frac{\alpha}{\rho}(1-\delta\rho L_g^2)(\tilde{\Xi}_{{\mathcal I}^c}-I)+\phi I\nonumber\\
&~~~+2\alpha\delta(I+\Xi_{\lambda,{\mathcal I}^c})J_{{\mathcal I}^c}J_{\mathcal I}^{\top}(J_{\mathcal I}J_{{\mathcal I}})^{-1}J_{\mathcal I}J_{{\mathcal I}^c}^{\top}(I+\Xi_{\lambda,{\mathcal I}^c})\nonumber\\
&\preceq\frac{\alpha}{\rho}(1-\delta\rho L_g^2)(\tilde{\Xi}_{{\mathcal I}^c}-I)+\phi I\nonumber\\
&~~~+2\alpha\delta(I+\Xi_{\lambda,{\mathcal I}^c})J_{{\mathcal I}^c}J_{{\mathcal I}^c}^{\top}(I+\Xi_{\lambda,{\mathcal I}^c})\nonumber\\
&\preceq\frac{\alpha}{\rho}(1-\delta\rho L_g^2)(\tilde{\Xi}_{{\mathcal I}^c}-I)+\phi I+8\alpha\delta B_g^2I,\label{equ54}
\end{align} where $A^{\top}(AA^{\top})^{-1}A\preceq I$ for a full row rank matrix $A$ has been applied in the second inequality.

If for all $i\in{\mathcal I}^c$ and $k\geq0$,
\begin{align}\label{equ55}
\tilde{\xi}_{i,k}\leq \pi^*=[\rho \max\limits_{i\in{\mathcal I}^c}\{g_i(x^*)\}/(\sqrt{C}d_0)+1]_+^2,
\end{align} then one can obtain that the sum on the right hand of (\ref{equ54}) is less than or equal to 0 by (\ref{equ14}) and (\ref{e16}), and then by Schur complement, $Q_2+\gamma I_m\preceq0$. Thus Theorem \ref{thm1} is proved.

What we left is to show (\ref{equ55}). Based on (\ref{equ9}) and $g_i(x^*)<0$ for $i\in{\mathcal I}^c$, one has that $\lambda_i^*=0$, $i\in{\mathcal I}^c$. For $k=0$, obviously, for all $i\in{\mathcal I}^c$, $\tilde{\xi}_{i,0}\leq\pi^*$, which indicates that (\ref{equ55}) holds for $k=0$, and then ${V}_{\delta,1}\leq V_{\delta,0}$. Therefore, invoking (\ref{equ48}) yields that for all $i\in{\mathcal I}^c$,
\begin{align*}
\tilde{\xi}_{i,1}\leq[\rho \max\limits_{i\in{\mathcal I}^c}\{g_i(x^*)\}/(\sqrt{C}d_0)+1]_+^2=\pi^*.
\end{align*} Subsequently, by the mathematical induction, (\ref{equ55}) can be proved. The proof is completed.
\hfill$\blacksquare$

\section{Example}\label{sec4}
In this section, an example motivated by applications in power systems \cite{tang2019semi} is presented to illustrate the feasibility of the discrete-time Aug-PGD (\ref{equ8}). Consider the following constrained optimization problem:
\begin{align}
&\min\limits_{p_i,q_i\in\mathbb{R}}~~~f(x)=\sum_{i=1}^n((p_i-p_{v,i})^2+q_i^2)\nonumber\\
&~~\text{s.t.}~~~~~~~~~~g_i(x)=p_1^2+q_i^2-S_i\leq0,\nonumber\\
&~~~~~~~~~~~~~~~~0\leq p_i\leq p_{v,i}, ~i\in[n],\label{equ58}
\end{align} where $x=(p_1,\ldots,p_n,q_1,\ldots,q_n)^{\top}$ and $p_{v,i},S_i$, $i\in[n]$ are constants. The problem (\ref{equ58}) along with an affine inequality constraint was considered in \cite{tang2019semi} but via a continuous-time dynamics Aug-PDGD. The affine inequality constraints can be regarded as special nonlinear constrains. Hence the algorithm Aug-PDG studied in this paper is applicable to the optimization problem (\ref{equ58}).

\begin{figure}[htbp]
\centering
\includegraphics[width=3in]{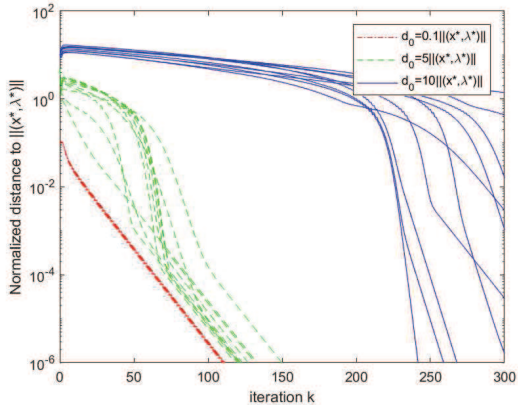}
\caption{Simulation of the relative distances to $\|(x^*,\lambda^*)\|$ with respect to iteration $k$.}\label{fig1}
\end{figure}

Let $n=10$,
\begin{align*}
S&=(S_1,\ldots,S_n)\\&=(2.7, 1.35, 2.7, 1.35, 2.025, 2.025, 2.7, 2.7, 1.35, 2.025)
 \end{align*} and $p_v=(p_{v,1},\ldots,p_{v,n})=4S$. Choose $\rho=0.1$ and $\alpha=0.1$. Three cases are simulated, where the initial point $(x_0,\lambda_0)$ is selected randomly such that the distance from the initial point $(x_0,\lambda_0)$ to the optimal point $(x^*,\lambda^*)$ (i.e., $d_0$) is $0.1\|(x^*,\lambda^*)\|$, $5\|(x^*,\lambda^*)\|$ and $10\|(x^*,\lambda^*)\|$, respectively. The curves of the normalized distance $\frac{\|(x_k-x^*,\lambda_k-\lambda^*)\|}{\|(x^*,\lambda^*)\|}$ with respect to the iteration $k$ are shown in Figure \ref{fig1}, where for each case, 10 instances of randomly selected initial points are considered. From Figure \ref{fig1}, it can be seen that the convergence rates are different for different $d_0$, and the distance $\|(x_k-x^*,\lambda_k-\lambda^*)\|$ linearly decays on the whole. Moreover, for each case,  the decreasing rate also changes as $(x_k,\lambda_k)$ approaches the optimal point. Specifically, the decreasing rates are small at the beginning and then become large when $(x_k,\lambda_k)$ goes to the optimal point. These observations support the semi-global linear convergence of the Aug-PDG, which is consistent with our theory analysis.

\section{Conclusion}\label{sec5}
In this paper, the linear convergence of an Aug-PDG in discrete-time for convex optimization with nonlinear inequality constraints has been investigated. Under some mild assumptions, the Aug-PDG has been proved to semi-globally converge at a linear rate, which depends on the distance from the initial point to the optimal point. Future research of interest may be to devise a method for solving the optimization problem under nonlinear constraints with a global linear convergence rate.

\bibliographystyle{IEEEtran}
\bibliography{MinMeng}

\end{document}